\documentclass[12pt]{article}

\usepackage{algorithmic}

\usepackage{amsthm,amsmath,amssymb}

\usepackage{stmaryrd}

\usepackage[colorlinks=true,citecolor=black,linkcolor=black,urlcolor=blue]{hyperref}

\usepackage{pstricks}
\newcommand{\showgrid}{}
\newcommand{\gridon}{\renewcommand{\showgrid}{\psset{subgriddiv=1,griddots=10,gridlabels=6pt}\psgrid}}
\gridon 

\usepackage{graphicx}

\usepackage{color}
\usepackage{multido}

\theoremstyle{plain}
\newtheorem{thm}{Theorem}

\newtheorem{cor}[thm]{Corollary}

\theoremstyle{definition}
\newtheorem{dfn}[thm]{Definition}
\newtheorem{ex}[thm]{Example}
\newtheorem{con}[thm]{Conjecture}

\theoremstyle{remark}

\newgray{hellgrau}{.89}

\def\bit{\begin{itemize}}
\def\eit{\end{itemize}}
\def\beq{\begin{equation}}
\def\eeq{\end{equation}}
\def\bpf{\begin{proof}}
\def\epf{\end{proof}}

\def\EM#1{{\em #1}}

\def\figref#1{Figure~\ref{#1}}

\def\Z{{\mathbb Z}}

\def\1{{\mathbf 1}}
\def\0{{\mathbf 0}}

\def\defeq{:=}

\def\of#1{\left(#1\right)}

\def\setof#1{\left\{#1\right\}}

\def\pas#1{\left(#1\right)}

\def\symm{\mathfrak S}

\def\inv{{\mathbf{inv}}}
\def\ASM{alternating sign matrix}
\def\ASMes{alternating sign matrices}
\def\DPP{descending plane partition}

\def\reell{real}
\def\virtuell{virtual}

\usepackage[applemac]{inputenc}

\title{\bf A simple bijection between permutation matrices and \DPP s
  without special parts.}

\author{Markus Fulmek\thanks{Research supported by the National Research Network ``Analytic
Combinatorics and Probabilistic Number Theory'', funded by the
Austrian Science Foundation.}\\
\small Faculty of Mathematics\\[-0.8ex]
\small University of Vienna\\[-0.8ex] 
\small Vienna, Austria\\
\small\tt Markus.Fulmek@univie.ac.at\\
}

\date{{2016--09--17}}

\begin{document}

\maketitle


\begin{abstract}
  We present a simple bijection between permutation matrices and \DPP s
  without special parts. This bijection is already mentioned in
  \cite{lalonde:2006} (without giving the details); it involves the inversion words of permutations
  and the (well--known) representation of \DPP s as families of non--intersecting lattice paths.
  
  (Taking a short detour, we will also exhibit how the (well--known) enumeration of \DPP s follows 
  easily from the evaluation of Andrew's determinant.)

\end{abstract}

\section{Introduction}
\label{sec:intro}

It is (nowadays) a well--known fact that the enumeration of \DPP s with parts not
exceding $m$ and of \ASMes\ of dimension $m$ gives the same number:
\begin{equation}
\label{eq:nof-ASMs}
\prod_{k=0}^{m-1}\frac{\pas{3\cdot k +1}!}{\pas{m + k}!}.
\end{equation}
On the other hand, it is apparently very hard to find a simple bijection between
\ASMes\ and \DPP s in general. Even the bijection between the much simpler special
cases of \EM{permutation matrices} (i.e., \ASMes\ without entries $-1$) and 
\DPP s with no \EM{special parts} (to be explained in a moment) led to complicated
constructions. In \cite{Ayyer:2010}, Ayyer presented an inductively constructed bijection
``which does not
relate the number of parts of the descending plane partition with the number of inversions of
the permutation as one might have expected from the conjecture
of Mills, Robbins and Rumsey'' (\cite[Conjecture 3]{MillsRobbinsRumsey:1983}, to be explained in a
moment). In \cite{Striker:2011}, Striker presented another bijection involving
monotone triangles as ``intermediate'' combinatorial objects, which maps
\DPP s with $k$ parts, all of which are non--special, to permutations with $k$ inversions.

In this note, we shall present a simple bijection which relies only on the
(obvious) representation of \DPP s as families of non--intersecting lattice
paths and on the (obvious) encoding of permutations by inversion words:
This bijection was mentioned in \cite{lalonde:2006} without giving the details.
It also  maps \DPP s with $k$ parts to permutations with $k$ inversions.

This note is organized as follows:

In section~\ref{sec:background}, we present the basic definitions and
background information. In section~\ref{sec:dpp-rep} we present the
interpretation of \DPP s as families of non--intersecting lattice paths.
In section~\ref{sec:nonspecial-bijection}, we present our simple
bijection.

\section{Background information}
\label{sec:background}

For reader's convenience, we recall some background information needed for our presentation.
\subsection{Descending plane partitions}
\label{sec:background-dpp}
Here is the definition of \DPP s (see \cite[Definition 4]{MillsRobbinsRumsey:1983}):
\begin{dfn}
\label{dfn:dpp}
A \EM{\DPP} is an array $\pi=\pas{a_{i,j}}$, $1\leq i<j<\infty$, of
\EM{positive integers} 
$$
\pi=\;\;
\begin{matrix}
a_{1,1} & a_{1,2} & a_{1,3} &         & \cdots       &        & \cdots      & a_{1,\mu_1} \\
        & a_{2,2} & a_{2,3} &         & \cdots       & \cdots & a_{2,\mu_2} & \\
        &         &         &         & \cdots       &        &             & \\
        &         &         &         & \cdots       &        &             & \\
        &         &         & a_{k,k} & \cdots &  a_{k,\mu_k} &             &
\end{matrix}
$$
such that
\begin{enumerate}
\item rows are \EM{weakly decreasing}, i.e., $a_{i,j}\geq a_{i,j+1}$ for all $i=1,\dots,k$ and $i\leq j < \mu_i$,
\item columns are \EM{strictly decreasing}, i.e., $a_{i,j}> a_{i+1,j}$ for all $i=1,\dots,k-1$ and $i<j\leq \mu_{i+1}$,
\item $a_{i,i}>\mu_i-i+1$ for all $i=1,\dots,k$,
\item $a_{i,i}\leq \mu_{i-1}-\pas{i-1}+1$ for all $i=2,\dots,k$.
\end{enumerate}

It is easy to see that these conditions imply
$$
\mu_1\geq\mu_2\geq\cdots\geq\mu_k\geq k.
$$
The \EM{parts} of a \DPP\ are the \EM{numbers} (with repetitions) that appear in the array.
The \EM{empty} array, which we denote by $\emptyset$, is explicitly allowed. 

A \DPP\ $\pi$ where no part is greater than $m$ (i.e., $\pi$ has at most $m-1$ rows)
is said to have dimension $m$.

We denote the $i$--th row of a \DPP\ by $r_i$.
The \EM{length} of $r_i$ is the number of parts it contains, which is $\mu_i-i+1$.
So we may rephrase the last two conditions as
\bit
\item[(A)] The first part of $r_i$ is greater than the length of $r_i$ for $i=1,\dots,k$,
\item[(B)] The first part of $r_i$ is less or equal than the length of the preceding row $r_{i-1}$  for $i=2,\dots,k$.
\eit

A part $a_{i,j}$ in a \DPP\ is called \EM{special} if it does not exceed the number of parts to
its left (in its row $r_i$), i.e.,
$$
a_{i,j}\leq j-i.
$$
\end{dfn}
\begin{ex}
A typical example is the array
$$
\begin{matrix}
6 & 6 & 6 & 4 & {\red\underline 2} \\
  & 5 & 3 & {\red\underline 2} & {\red\underline 1} \\
  &   & 2 &   & 
\end{matrix}
$$
with $3$ rows and $10$ parts (written in descending order)
$$
6,6,6,5,4,3,2,{\red\underline 2,\underline 2,\underline 1},
$$
three of which are  special parts (indicated as underlined numbers;
note that the $2$ in the last row is \EM{not}
a special part):
$$
{\red \underline 2,\underline 2,\underline 1}.
$$
(This is example $D_0$ in 
\cite{Lalonde2003311}, see \cite[Fig.~1]{Lalonde2003311}.)

\end{ex}
\subsection{Alternating sign matrices}
Here is the definition of \ASMes\ (see \cite[Definition 1]{MillsRobbinsRumsey:1983}):
\begin{dfn}
An \EM{\ASM} of dimension $m$ is an $m\times m$ square matrix which satisfies
\bit
\item all entries are $1$, $-1$ or $0$,
\item every row and column has sum $1$,
\item in every row and column the nonzero entries alternate in sign.
\eit

Suppose that $M=\pas{a_{i,j}}_1^m$ is an alternating sign matrix of dimension $m$.
Then the number of \EM{inversions} in $M$ is defined to be
\begin{equation}
\label{eq:inversions}
\sum_{\substack{1\leq i<k\leq n\\1\leq l<j\leq n}}a_{i,j}\cdot a_{k,l}.
\end{equation}
(See \cite[p.\ 344]{MillsRobbinsRumsey:1983}.)
\end{dfn}

\begin{ex}
The following matrix is an example of an \ASM\ of dimension $5$:
$$
\begin{pmatrix}
0 & 1 & 0 & 0 & 0 \\
0 & 0 & 1 & 0 & 0 \\
1 & -1 & 0 & 1 & 0 \\
0 & 1 & 0 & -1 & 1 \\
0 & 0 & 0 & 1 & 0.
\end{pmatrix}
$$
\end{ex}

\subsection{The Mills--Robbins--Rumsey conjecture}
Here is  the Conjecture of Mills, Robbins and Rumsey \cite[Conjecture 3]{MillsRobbinsRumsey:1983}:
\begin{con}
\label{con:MRR}
Suppose that $m,k,n,p$ are nonnegative integers, $1\leq k\leq m$.
Let $\mathcal{A}\of{m,k,n,p}$ be the set of \ASMes\ such that
\begin{enumerate}
\item the size of the matrix is $m\times m$,
\item the $1$ in the top row occurs in position $k$,
\item the number of $-1$'s in the matrix is $n$,
\item the number of inversions in the matrix is $p$.
\end{enumerate}
On the other hand, let $\mathcal{D}\of{m,k,n,p}$ be the set of \DPP s such that
\begin{enumerate}
\item no part exceeds $m$,
\item there are exactly $k-1$ parts equal to $m$,
\item there are exactly $n$ special parts,
\item there are a total of $p$ parts.
\end{enumerate}
Then $\mathcal{A}\of{m,k,n,p}$ and $\mathcal{D}\of{m,k,n,p}$ have the same cardinality.
\end{con}

\subsection{Permutation matrices and inversions}
Let $\sigma\in\symm_m$ be a permutation of the first $m$ natural numbers $\setof{1,2,\dots,m}$.
Recall that an \EM{inversion} of $\sigma$ is a pair $\pas{i,j}$ such that
$i<j$ but $\sigma\of i > \sigma\of j$. (For the number $\inv\of\sigma$ of all inversions of
$\sigma$ we have $0\leq\inv\of\sigma\leq \frac{m\cdot\pas{m+1}}2$.)

We may assign to $\sigma$ its \EM{inversion word} $\pas{a_1, a_2,\dots, a_{m-1}}$, where
$a_k$ is the number of inversions $\pas{i,j}$ with $\sigma\of j=k$, $k = 1,2,\dots, m-1$.
Clearly we have $0\leq a_k\leq m-k$ and $a_1+a_2+\cdots+a_{m-1} = \inv\of\sigma$.

Considering the \EM{permutation word}
$$
\pas{\sigma\of 1,\sigma\of 2,\dots,\sigma\of m},
$$
of $\sigma$, the inversion word's $k$--th entry $a_k$ is simply the number of elements
\EM{to the left} of $k$ (in the permutation word) which are \EM{greater} than $k$, and
it is easy to see that \EM{every} word $\pas{b_1,b_2,\dots,b_{m-1}}$ with $0\leq b_k\leq m-k$
determines a \EM{unique} permutation: Inversion words are, in this sense, just another
``encoding'' for permutations.

A permutation $\sigma\in\symm_n$ can be represented by an $m\times m$--matrix $M$ with
entries $0$ or $1$, namely
$$
M_{i,j} = \delta_{\sigma\of j,i}
$$
(where $\delta_{x,y}$ denotes Kronecker's delta: $\delta_{x,y} = 1$ if $x=y$,
$\delta_{x,y} = 0$ if $x\neq y$). We call this matrix the \EM{permutation matrix}
of $\sigma$: Clearly, it contains precisely one entry $1$ in every row and column.
\begin{ex}
Let $\sigma\in\symm_6$ be the permutation with permutation word 
$$
6 3 1 4 2 5.
$$
The corresponding permutation matrix is
$$
\begin{pmatrix}
0 & 0 & 1 & 0 & 0 & 0 \\
0 & 0 & 0 & 0 & 1 & 0 \\
0 & 1 & 0 & 0 & 0 & 0 \\
0 & 0 & 0 & 1 & 0 & 0 \\
0 & 0 & 0 & 0 & 0 & 1 \\
1 & 0 & 0 & 0 & 0 & 0
\end{pmatrix},
$$
and the corresponding inversion word is
$$
\pas{2,3,1,1,1}.
$$
\end{ex}

Note that every \EM{permutation matrix} is also an \ASM\ (which does not contain entries
$-1$), and the definition of inversions \eqref{eq:inversions} for \ASMes\ is a generalization
of the number of inversions of a permutation.

\section{Representation of \DPP s as lattice paths}
\label{sec:dpp-rep}
\begin{figure}
\caption{Illustration}\label{fig:lalonde-example}
Consider the following \DPP\ $\pi$ of dimension $6$:
$$
\begin{matrix}
6 & 6 & 6 & 4 & {\red\underline 2}\\
  & 5 & 3 & {\red\underline 2} & {\red\underline 1}\\
  &   & 2 &  & \\
\end{matrix}
$$
The picture shows the visualization of $\pi$ as a family of $5=6-1$
non--intersecting lattice paths: 

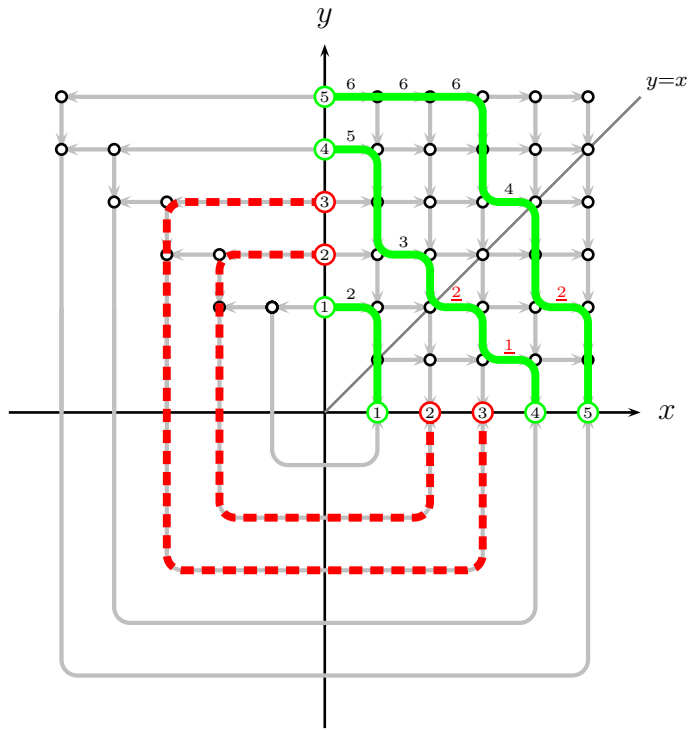
\begin{center}

\input save-graphics/DPP-alternativ.tex
\end{center}
Here, the starting point $S_i$ corresponds to the point numbered $i$
on the vertical axis, and the ending point $E_i$ corresponds to the
point numbered $i$ on the horizontal axis. The $3$ \reell\ paths corresponding
to the $3$ rows of $\pi$ are shown as solid lines coloured green, and the $2$ remaining \virtuell\
paths (connecting $S_2$ to $E_2$ and $S_3$ to $E_3$) are shown as dashed lines coloured red.
Note that the special parts of $\pi$ correspond to the horizontal steps
in the ``special range'' below the main diagonal $y=x$. 
\end{figure}

We shall use the well--known encoding of (shifted) tableaux as non--intersec\-ting
lattice paths, with small additions/modifications; i.e., we shall encode \DPP\
$\pi=\pas{a_{i,j}}$ of dimension $m$ as a set
of $m-1$ non--intersecting lattice paths in a particular sub--lattice
${\mathcal L}_m$ of $\Z^2$. This encoding was already considered in \cite{Lalonde2003311},
but for reader's convenience, we shall describe its details below. Certainly, it
will be helpful to look at an illustrative
example; see \figref{fig:lalonde-example}.

The starting points of these lattice paths are the points
$$
\mathcal{S}=\setof{S_1\defeq\pas{0,2},S_2\defeq\pas{0,3},\dots,S_{m-1}\defeq\pas{0,m}}\subset\Z^2
$$
on the vertical axis.

The ending points of these lattice paths are the points 
$$
\mathcal{E}=\setof{E_1\defeq\pas{1,0},E_2\defeq\pas{2,0},\dots,E_{m-1}\defeq\pas{m-1,0}}\subset\Z^2
$$
on the horizontal axis.

Moreover, the following \EM{rectangle of points} also belongs to ${\mathcal L}_m$:
$$
\mathcal{R}=\setof{1,2\dots,m-1}\times\setof{1,2,\dots,m}
\subset\Z^2.
$$

If two points $\mathbf p_1=\pas{x_1,y_1}$, $\mathbf p_2=\pas{x_2,y_2}$
in the union of
$$
\mathcal{S}\cup\mathcal{E}\cup\mathcal{R},
$$
which are \EM{not} both starting points \EM{and not} both end points,
are ``horizontally or vertically adjacent'', then there is a (directed)
arc ``from left to right or from top to bottom'' between them, i.e.,
\bit
\item if ${\mathbf{p}}_2=\pas{x_2,y_2} = \pas{x_1+1,y_1}$, then there is an arc
from $\mathbf p_1$ to $\mathbf p_2$,
\item if ${\mathbf p}_2=\pas{x_2,y_2} = \pas{x_1,y_1-1}$, then there is an arc
from $\mathbf p_1$ to $\mathbf p_2$.
\eit

We shall call the directed graph consisting of these points and arcs
the \EM{\reell\ part} of our lattice ${\mathcal L}_m$, and we use it in the
obvious way to encode the rows of $\pi$ as a family of \EM{\reell\ non--intersecting
lattice paths}: Row $i$
of $\pi$ directly corresponds to a path in (the \reell\ part of)
the lattice ${\mathcal L}_m$ by interpreting the $j$--th \EM{part}
$a_{i,j+i-1}$ of this row
as the \EM{height} of the $j$--th horizontal step of the corresponding path
(see \figref{fig:lalonde-example}): Note that this path has starting
point $S_{a_{i,i}-1}$, which is \EM{not necessarily equal} to $S_{m-i}$.

We define the \EM{length} of some \reell\ path to be the number of
its horizontal steps: So the length of the path is equal to the length
of the row it encodes.

In general, this ``interpreting rows as lattice paths'' will not give
$m-1$ paths: The missing (\virtuell) paths shall use the \EM{\virtuell\ part} of our
lattice ${\mathcal L}_m$, which contains the starting and ending points
$\mathcal{S}\cup\mathcal{E}$ together with the union of two
``diagonals'', namely
\bit
\item $\mathcal{D}_1=
\setof{\pas{-x,x}:\;2 \leq x \leq m-1}$,
\item $\mathcal{D}_2=
\setof{\pas{-x,x+1}:\;1 \leq x \leq m-1}$.
\eit
The (directed) arcs in this \virtuell\ part of our lattice ${\mathcal L}_m$
are
\bit
\item from $S_i=\pas{0,i+1}$ to $\pas{-i,i+1}$ for $i=1,2,\dots,m-1$,
\item from $\pas{-i,i+1}$ to $\pas{-i-1,i+1}$ for $i=1,2,\dots,m-2$,
\item from $\pas{-i-1,i+1}$ to $E_{i+1}=\pas{i+1,0}$ for $i=1,2,\dots,m-2$,
\item from $\pas{-1,2}$ to $E_1=\pas{1,0}$.
\eit
Note that there are \EM{unique} \virtuell\ lattice paths connecting the points
\bit
\item $S_i$ and $E_i$, for $i=1,2,\dots,m-1$,
\item $S_i$ and $E_{i+1}$ for $i=1,2,\dots,m-2$.
\eit

If there are starting points and ending points which are not connected by
\reell\ paths, then we shall connect them with non--intersecting \virtuell\ 
paths. 
In fact, this is always possible in a \EM{unique} way:

Condition 3, rephrased as condition (A) in section~\ref{sec:background-dpp}, in the definition of a \DPP\ states that
$a_{k,k}$ is strictly greater than the length of the \reell\ path starting at height $a_{k,k}$, so a \reell\ path starting at point $S_{i}=\pas{0,i+1}$ has length
$j \leq i$ and therefore must end in some point $E_{j}=\pas{j,0}$ with $j\leq i$.

Condition 4, rephrased as condition (B) in section~\ref{sec:background-dpp}, in the definition of a \DPP\ states that
$a_{k,k}$ is less or equal than the length of the path starting at height $a_{k-1,k-1}$, so
a \reell\ path starting in $S_{i}=\pas{0,i+1}$ and ending in some point
$E_{j}=\pas{j,0}$ with $j\leq i$ implies that the \reell\ path \EM{below}
(if any) starts in some starting point $S_{\ell}$ with $\ell<j$.

So, if there are $d$ \reell\ lattice paths connecting the points
$S_{s_i}$ and $E_{e_i}$ for $i=1,2,\dots d$ in their ``natural order''
(i.e., $s_d> s_{d-1}>\cdots> s_1$ and  $e_d> e_{d-1}>\cdots> e_1$) then there holds
$$
s_d\geq e_d > s_{d-1}\geq e_{d-1} > \cdots > s_2\geq e_2 > s_1 \geq e_1.
$$

Now assume that all starting and ending points $S_i$, $E_i$, $a< i\leq m-1$,
are properly connected by \reell\ or \virtuell\ paths as described above,
and that the next row $r$ of the \DPP\ (not yet encoded as a \reell\ path)
starts with part $b\leq a$ (at the beginning, this assumption is fulfilled for
$a=m-1$). If $b<a$, then for $i=b+1,b+2,\dots a$ we connect starting point
$S_i$ with ending point $E_i$ by a \virtuell\ path: Clearly, there is one
and only one way to achieve this. Let $c\leq b$ be the length of row $r$
and connect $S_b$ and $E_c$ with the
\reell\ lattice path corresponding to row $r$. If $c<b$, then for $i=c,c+1,\dots,b-1$ we connect $S_i$ with $E_{i+1}$
by a \virtuell\ path: Clearly, there is one and only one way to achieve this.
Observe that by now we achieved a proper connection of all starting and ending
points $S_i$, $E_i$ for $c\leq i\leq m-1$.

Repeating this step for every row $r$ of the given \DPP\ might leave $q$
pairs of starting points and ending points $\pas{S_i, E_i}$, $i=1,2,\cdots q<m$,
which are not yet connected properly: If $q>0$, then for $i=1,2,\cdots q$,
we connect starting point $S_i$ with ending point $E_i$ by a \virtuell\ path
in the only possible way.

It is easy to see that the $m-1$ lattice paths thus constructed are
\EM{non--intersecting}, so we obtained an encoding
of a \DPP\ of dimension $m$ as a family of $m$ non--intersecting lattice paths
in ${\mathcal L}_m$.

We call starting (or ending) points \reell\ (\virtuell) if they belong to a \reell\
(\virtuell) path.

Now we show that \EM{every} family $F$ of non--intersecting lattice paths
connecting the starting points and the ending points in our lattice
${\mathcal L}_m$ determines a \EM{unique} \DPP: Arranging the heights $h_{i,j+i-1}$ of the $j$--th horizontal step
of the $i$--th \reell\ path (counted from above) as follows \dots
$$
\begin{matrix}
h_{1,1} & a_{1,2} & h_{1,3} &         & \cdots       &        & \cdots      & h_{1,\mu_1} \\
        & h_{2,2} & h_{2,3} &         & \cdots       & \cdots & h_{2,\mu_2} & \\
        &         &         &         & \cdots       &        &             & \\
        &         &         &         & \cdots       &        &             & \\
        &         &         & h_{l,l} & \cdots &  h_{l,\mu_l} &             &
\end{matrix}
$$
\dots gives a (unique) array of positive integers (maybe empty) which fulfils
conditions (1) and (2)
in Definition~\ref{dfn:dpp}, since the \reell\ paths are non--intersecting. In order
to show that this array is indeed a \DPP, we must check conditions (3) and (4) (rephrased
as (A) and (B) in section~\ref{sec:background-dpp}).

Assume that $F$ contains a \reell\ lattice path $p$ connecting $S_a$ with 
$E_b$. Consider the set of (all) lattice paths ``above $p$'', i.e., starting in points $S_i$, for
$a<i<m$.
Let $\ell$ be the number of \reell\ paths in this set. Then this set contains
$\pas{m-1-a} - \ell$ \virtuell\ paths, which, of course,  must have $\pas{m-1-a} - \ell$ \virtuell\
ending points:
By construction, for such \virtuell\ ending point $E_j$ we must have $j>a$. If $b >a$, then
$\ell+1$ of these possible ending points would, in fact, be \reell\ ending points
(since the paths are non--intersecting),
leaving only $\pas{m-1-a} - \ell - 1$ possible \virtuell\ ending points, which simply is
not enough. Hence there must hold $b\leq a$, which is equivalent to condition (3)
(or rephrased condition (A)) from Definition~\ref{dfn:dpp}.

Now consider the set of (all) lattice paths ending in points $E_i$, for $b<i<m$.
Again, let $\ell$ be the number of \reell\ paths in this set, so there must be 
$\pas{m-1-b} - \ell$ \virtuell\ paths in this set, which, of course,  must have
$\pas{m-1-b} - \ell$ \virtuell\ starting points: By construction, for every such \virtuell\ 
starting point $S_j$ we must have $j\geq b$. Since $a\geq b$, among these possible
starting points there are at least $\ell+1$ \reell\ points, leaving at most
$\pas{m-1}-\pas{b-1}-\pas{\ell+1}$ possible points, which is precisely the required
number. Therefore, for every \reell\ path with starting point $S_c$ with $c<a$ there
must hold $c<b$, which is equivalent to condition 4
(or rephrased condition (B)) from Definition~\ref{dfn:dpp}.

So we established the
bijection between \DPP s of dimension
$m$ and families of $\pas{m-1}$ non--intersecting lattice paths in our lattice
${\mathcal L}_m$.

\subsection{Detour: A determinantal formula}
From this representation, we immediately obtain the following determinantal
expression for the number of \DPP s:
\begin{cor}
The number of \DPP s of dimension $m$ is 
\begin{equation}
\label{eq:dpp-det}
\det\pas{\binom{i+j-1}{j-1}+\delta_{i,j}-\delta_{i+1,j}}_{1\leq i,j\leq m-1}.
\end{equation}
\end{cor}
\begin{proof}
Note that 
$$
a_{i,j}\defeq\binom{i+j-1}{j-1}+\delta_{i,j}+\delta_{i+1,j}
$$
is precisely the number of lattice paths in ${\mathcal L}_m$ connecting
starting point $S_i$ and ending point $E_j$.

The straightforward application
of the \EM{Lindström--Gessel--Viennot Theorem}\footnote{Using Krattenthaler's
\cite[footnote 10 on page 76]{Krattenthaler:2005}
name for this well--known result.} (see  \cite{GesselViennot:1985,GesselViennot:1989})
shows that in the expansion of the determinant $\det\of{a_{i,j}}$, all terms corresponding
to \EM{intersecting} families of $\pas{m-1}$ lattice paths cancel;
and each \EM{non--intersecting} family
connecting 
$S_i$ to $E_{\sigma\of i}$, $i=1,\dots,m-1$, will be counted with the \EM{sign}
of the permutation $\sigma$.
Unfortunately, due to the particular construction of our lattice, there are
permutations with a negative sign which ``survive'' the Lindström--Gessel--Viennot--cancellation
of intersecting lattice paths, and the corresponding families
of lattice paths will be subtracted from instead of added to the number we want to
determine.

But it is easy to see
that the sign of the permutation
$\sigma$ is precisely $\pas{-1}^k$, where $k$ is the number of \virtuell\ paths connecting
some starting point $S_i$ with an ending point $E_{i+1}$. So giving all such paths
weight $\pas{-1}$, i.e., considering the corresponding ``weighted'' number of
lattice paths
$$
b_{i,j}\defeq\binom{i+j-1}{j-1}+\delta_{i,j}-\delta_{i+1,j}
$$
instead of $a_{i,j}$, will cancel out all the negative signs; and thus
the determinant \eqref{eq:dpp-det} provides the correct number of
\DPP s of dimension $m$.
\end{proof}

\subsection{Detour, continued: Andrew's determinant}
Note that this determinant is closely related to the famous determinant considered by
Andrews \cite{Andrews:1979}:
$$
a_{n}\of x\defeq\det\pas{\binom{x+i+j}{j}+\delta_{i,j}}_{0\leq i,j\leq n-1}.
$$
The analogous generalization of the determinant in \eqref{eq:dpp-det} would be
$$
d_n\of x\defeq\det\pas{\binom{x+i+j}{j}+\delta_{i,j}-\delta_{i+1,j}}_{0\leq i,j\leq n-1}.
$$
Obviously, $d_{m-1}\of 1$ is precisely the determinant in \eqref{eq:dpp-det}.

But since
$$
\pas{\binom{x+i+j}{j}+\delta_{i,j}}
-
\pas{\binom{x+i+\pas{j-1}}{j-1}+\delta_{i,j-1}}
$$
equals
$$
\binom{\pas{x-1}+i+j}{j}+\delta_{i,j}-\delta_{i+1,j}
$$
for all $j\geq 1$, we see that
$$
d_n\of x = a_n\of{x+1}\text{ for all } n>0.
$$
So we obtained that the number of \DPP s of dimension $m$ is equal to $a_{m-1}\of 2$.

Andrews (\cite[Theorem 8]{Andrews:1979}; see \cite{AndrewsStanton:1998} for a short proof) showed that 
$$
a_m\of x =
\prod_{k=0}^{m-1}\Delta_k\of x,
$$
where $\Delta_0\of x\equiv 2$ and for all $j>0$
\begin{align*}
\Delta_{2j}\of x   &=
\frac{
\pas{x+2j+2}_j\pas{x/2+2j+3/2}_{j-1}
}{
\pas{j}_j\pas{x/2+j+3/2}_{j-1}
} , \\
\Delta_{2j-1}\of x &= 
\frac{
\pas{x+2j}_{j-1}\pas{x/2+2j+1/2}_{j}
}{
\pas{j}_j\pas{x/2+j+1/2}_{j-1}
}.
\end{align*}
(Here, we used Pochhammer's symbol: $\pas{x_j}=x\cdot\pas{x+1}\cdots\pas{x+j-1}$.)

Note that
$$
\Delta_{k-1}\of 2 = \frac{\pas{3 k + 1}!}{\pas{2 k + 1}!\pas{k + 1}_{k}}
$$
for all $k\geq 1$, and 
$$
\frac{\pas{3\cdot 0 + 1}!}{\pas{2\cdot 0 + 1}!\pas{1}_{0}} = 1,
$$ 
whence 
we obtain the following expression for the number of
\DPP s of dimension $m$:
$$
a_{m-1}\of 2 = \prod_{k=0}^{m-1}\frac{\pas{3k+1}!}{\pas{2k+1}!\pas{k+1}_k}.
$$
It is easy to see that this, in fact, is equal to \eqref{eq:nof-ASMs}.

\section{The bijection between \DPP s without special parts and permutations}
\label{sec:nonspecial-bijection}

Now we shall present the promised bijection between \DPP s and inversion words (which
are in bijection with permutations and with permutation matrices, as outlined in section~\ref{sec:background}):
 
If we are given a \DPP\ $\pi$ of dimension $m$ without special parts, we can easily derive from it the
\EM{inversion word} $\pas{a_1,a_2,\dots,a_{m-1}}$ of a permutation $\sigma\in\symm_m$: Simply
set
$$
a_i\defeq {\text{number of parts }\pas{m - i + 1}\text{ in }\pi}.
$$
Looking at the representation of $\pi$ as a family of non--intersecting lattice paths
in ${\mathcal L}_m$, it is easy to see that the number of (non--special) horizontal steps at
height $h$ cannot exceed $h-1$, whence we have $0\leq a_i\leq m-i$ for all $i$. So the
word $\pas{a_1,a_2,\dots,a_{m-1}}$ is indeed an inversion word which encodes some
unique permutation $\sigma\in\symm_m$. 

The inverse mapping
is also quite simple: If we have an inversion word $\pas{a_1,a_2,\dots,a_{m-1}}$, we
start with the empty \DPP\ and insert successively (i.e., for , $i=1,2,\dots,m-1$) $a_i$
parts $\pas{m-i+1}$, into the rows of a ``growing'' \DPP, subject to the simple rule, that
we never start a new row  that would violate condition (4) (rephrased as (B)) from
Definition~\ref{dfn:dpp}: Table~\ref{table:psi}
gives the corresponding algorithm in ``pseudo--code'' notation.
\begin{table}
\caption{Pseudo--code for the mapping from inversion words to \DPP s.}
\label{table:psi}
\begin{center}
\begin{algorithmic}
\COMMENT{Construct $m$-DPP from inversion word $w=\pas{w_1,\dots,w_{m-1}}$.}

\COMMENT{Input: $w=\pas{w_1,\dots,w_{m-1}}$ ($0\leq w_i\leq m-i$).}

\COMMENT{Output: DPP of dimension $m$ (as array of rows).}
\STATE $i\gets 0$ \COMMENT{$i$: current index of inversion word $w$.}
\STATE $r\gets 1$ \COMMENT{$r$: index of the (yet empty) row to be filled.}

\REPEAT
\STATE $i\gets i+1$
\COMMENT{Find next index $i$ such that $w_i>0$:}

\WHILE{$i < m$ \AND $w_i = 0$}
	 \STATE $i\gets i+1$
\ENDWHILE

\IF{$i\geq m$} \STATE {\bf break} \COMMENT{Jump out of loop} \ENDIF
\STATE $e\gets m - i + 1$ \COMMENT{We have to insert $w_i$ steps at height $e$:}
\STATE $l\gets\pas{\text{length of row }r}$
\COMMENT{Do not start a new row if $e>l$!}
\IF{$l\leq e$}
\STATE $r\gets r+1$ \COMMENT{Start a new row: Room for $m-i$ steps!}
\STATE $\pas{\text{start new row $r$ and insert $w_i$ entries $e$ into it}}$
\ELSE
\STATE $a = \min\of{e-l,w_i}$
\STATE $\pas{\text{append $a$ entries $e$ to row $r$}}$
\STATE $a\gets w_i-a$
\IF{$a>0$}
\STATE $r\gets r+1$ \COMMENT{Start a new row:  Room for $m-i$ steps!}
\STATE $\pas{\text{start new row and insert $a$ entries $e$ into it}}$
\ENDIF
\ENDIF
\UNTIL $i\geq m$.
\RETURN $\pas{\text{array of rows thus constructed}}.$
\end{algorithmic}
\end{center}
\end{table}
The correctness of this ``insertion of horizontal steps corresponding to the
parts of some \DPP\ \EM{without} special parts''
is easily seen by observing 
\bit
\item that the real path
``currently under construction'' 
must reach the  line $y=x$ \EM{before} the ``next'' (i.e., ``lower'') real path may start
(according to condition (4), rephrased as (B), from
Definition~\ref{dfn:dpp} for \DPP s);
\item and if some path reached the line $y=x$, then it has a unique continuation (by
vertical steps only, since there are no special parts) and can thus be ``finished''.
\eit

Clearly, this bijection is in line with the Mills--Robbins--Rumsey--Conjec\-tu\-re
(given here as Conjecture~\ref{con:MRR}): The number of parts
of the \DPP\ equals the number of inversions of the permutation, and if the number of parts
of the \DPP\ which are equal to $m$ is $k-1$, then the position of the $1$ in the first
row of the permutation matrix is $k$.

We conclude this presentation with an illustrating example:
\begin{ex}
Consider the inversion word
$$
a=\pas{0,0,2,3,1,1,1}.
$$
We shall illustrate the algorithm by showing the ``successively growing'' family
of (real) non--intersecting lattice paths corresponding to a \DPP\ of dimension $8$
without special parts. We start with the empty lattice. Since $a_1=a_2=0$, the first
horizontal steps to be inserted are $a_3=2$ steps at height $8-3+1=6$:
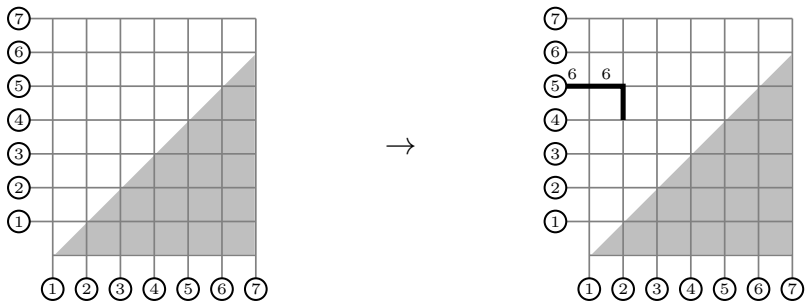
\begin{center}
\psset{unit=0.45cm}
\pspicture(-1,-1)(8,9)
\pspolygon[linecolor=white,fillstyle=solid,fillcolor=lightgray](1,1)(7,1)(7,7)
\psset{linewidth=0.05,linecolor=gray}
\multido{\nx=1+1}{7}{\psline(\nx,0)(\nx,8)}
\multido{\ny=2+1}{7}{\psline(0,\ny)(7,\ny)}
\psline(1,1)(7,1)

\psset{fillstyle=none,linecolor=black,linewidth=0.15}
\psset{linewidth=0.06,linecolor=black,fillstyle=solid,fillcolor=white}
\multido{\nx=1+1}{7}{\pscircle(\nx,0){0.35}\rput(\nx,0){{\tiny \nx}}}
\multido{\ny=2+1}{7}{\pscircle(0,\ny){0.35}}
\multido{\ny=1+1}{7}{\rput(0,\ny){\rput(0,1){\tiny\ny}}}

\endpspicture

\hfil
\hbox{\raisebox{2.25cm}{$\to$}}
\hfil
\psset{unit=0.45cm}
\pspicture(-1,-1)(8,9)
\pspolygon[linecolor=white,fillstyle=solid,fillcolor=lightgray](1,1)(7,1)(7,7)
\psset{linewidth=0.05,linecolor=gray}
\multido{\nx=1+1}{7}{\psline(\nx,0)(\nx,8)}
\multido{\ny=2+1}{7}{\psline(0,\ny)(7,\ny)}
\psline(1,1)(7,1)

\psset{fillstyle=none,linecolor=black,linewidth=0.15}
\psline(0,6)(2,6)(2,5)
\rput(0.5,6.35){{\tiny 6}}
\rput(1.5,6.35){{\tiny 6}}
\psset{linewidth=0.06,linecolor=black,fillstyle=solid,fillcolor=white}
\multido{\nx=1+1}{7}{\pscircle(\nx,0){0.35}\rput(\nx,0){{\tiny \nx}}}
\multido{\ny=2+1}{7}{\pscircle(0,\ny){0.35}}
\multido{\ny=1+1}{7}{\rput(0,\ny){\rput(0,1){\tiny\ny}}}

\endpspicture

\end{center}
Now we have to insert $a_4=3$ steps at height $8-4+1=5$: Since we did not reach the
line $y=x$ yet, we must not start a new path, but append these steps to the current path
--- by doing this, we reach the line $y=x$ and are thus able to ``finish'' this path
(since we must not insert horizontal steps \EM{below} $y=x$).
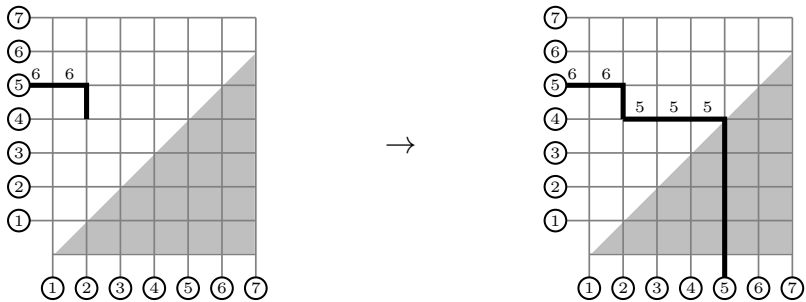
\begin{center}

\hfil
\hbox{\raisebox{2.25cm}{$\to$}}
\hfil
\psset{unit=0.45cm}
\pspicture(-1,-1)(8,9)
\pspolygon[linecolor=white,fillstyle=solid,fillcolor=lightgray](1,1)(7,1)(7,7)
\psset{linewidth=0.05,linecolor=gray}
\multido{\nx=1+1}{7}{\psline(\nx,0)(\nx,8)}
\multido{\ny=2+1}{7}{\psline(0,\ny)(7,\ny)}
\psline(1,1)(7,1)

\psset{fillstyle=none,linecolor=black,linewidth=0.15}
\psline(0,6)(2,6)(2,5)
\rput(0.5,6.35){{\tiny 6}}
\rput(1.5,6.35){{\tiny 6}}
\psline(2,5)(5,5)(5,4)(5,0)
\rput(2.5,5.35){{\tiny 5}}
\rput(3.5,5.35){{\tiny 5}}
\rput(4.5,5.35){{\tiny 5}}
\psset{linewidth=0.06,linecolor=black,fillstyle=solid,fillcolor=white}
\multido{\nx=1+1}{7}{\pscircle(\nx,0){0.35}\rput(\nx,0){{\tiny \nx}}}
\multido{\ny=2+1}{7}{\pscircle(0,\ny){0.35}}
\multido{\ny=1+1}{7}{\rput(0,\ny){\rput(0,1){\tiny\ny}}}

\endpspicture

\end{center}
Now we have to insert $a_5=1$ step at height $8-5+1=4$: Since the ``preceding''
path is finished, we start a new one.
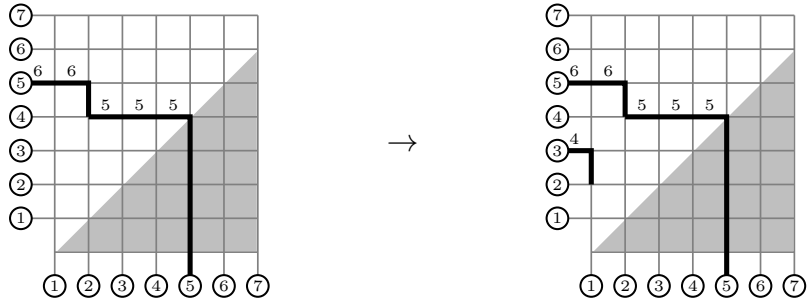
\begin{center}

\hfil
\hbox{\raisebox{2.25cm}{$\to$}}
\hfil
\psset{unit=0.45cm}
\pspicture(-1,-1)(8,9)
\pspolygon[linecolor=white,fillstyle=solid,fillcolor=lightgray](1,1)(7,1)(7,7)
\psset{linewidth=0.05,linecolor=gray}
\multido{\nx=1+1}{7}{\psline(\nx,0)(\nx,8)}
\multido{\ny=2+1}{7}{\psline(0,\ny)(7,\ny)}
\psline(1,1)(7,1)

\psset{fillstyle=none,linecolor=black,linewidth=0.15}
\psline(0,6)(2,6)(2,5)
\rput(0.5,6.35){{\tiny 6}}
\rput(1.5,6.35){{\tiny 6}}
\psline(2,5)(5,5)(5,4)(5,0)
\rput(2.5,5.35){{\tiny 5}}
\rput(3.5,5.35){{\tiny 5}}
\rput(4.5,5.35){{\tiny 5}}
\psline(0,4)(1,4)(1,3)
\rput(0.5,4.35){{\tiny 4}}
\psset{linewidth=0.06,linecolor=black,fillstyle=solid,fillcolor=white}
\multido{\nx=1+1}{7}{\pscircle(\nx,0){0.35}\rput(\nx,0){{\tiny \nx}}}
\multido{\ny=2+1}{7}{\pscircle(0,\ny){0.35}}
\multido{\ny=1+1}{7}{\rput(0,\ny){\rput(0,1){\tiny\ny}}}

\endpspicture

\end{center}
Now we have to insert $a_6=1$ step at height $8-5+1=3$: Since the ``preceding''
path is not yet finished, we append this step to it: We see that this path
will reach the line $y=x$ in the next iteration, so we may ``finish'' it already.
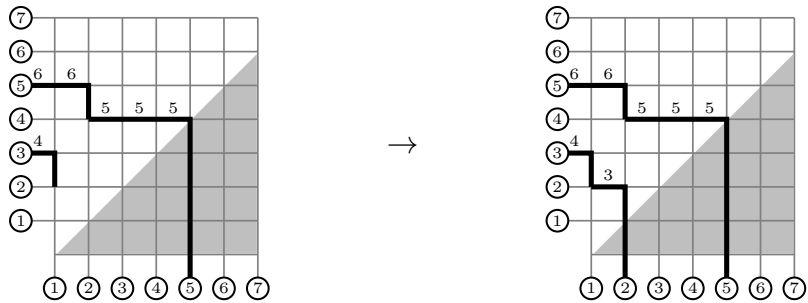
\begin{center}

\hfil
\hbox{\raisebox{2.25cm}{$\to$}}
\hfil
\psset{unit=0.45cm}
\pspicture(-1,-1)(8,9)
\pspolygon[linecolor=white,fillstyle=solid,fillcolor=lightgray](1,1)(7,1)(7,7)
\psset{linewidth=0.05,linecolor=gray}
\multido{\nx=1+1}{7}{\psline(\nx,0)(\nx,8)}
\multido{\ny=2+1}{7}{\psline(0,\ny)(7,\ny)}
\psline(1,1)(7,1)

\psset{fillstyle=none,linecolor=black,linewidth=0.15}
\psline(0,6)(2,6)(2,5)
\rput(0.5,6.35){{\tiny 6}}
\rput(1.5,6.35){{\tiny 6}}
\psline(2,5)(5,5)(5,4)(5,0)
\rput(2.5,5.35){{\tiny 5}}
\rput(3.5,5.35){{\tiny 5}}
\rput(4.5,5.35){{\tiny 5}}
\psline(0,4)(1,4)(1,3)
\rput(0.5,4.35){{\tiny 4}}
\psline(1,3)(2,3)(2,0)
\rput(1.5,3.35){{\tiny 3}}
\psset{linewidth=0.06,linecolor=black,fillstyle=solid,fillcolor=white}
\multido{\nx=1+1}{7}{\pscircle(\nx,0){0.35}\rput(\nx,0){{\tiny \nx}}}
\multido{\ny=2+1}{7}{\pscircle(0,\ny){0.35}}
\multido{\ny=1+1}{7}{\rput(0,\ny){\rput(0,1){\tiny\ny}}}

\endpspicture

\end{center}
Finally, we have to insert $a_7=1$ step at height $8-7+1=2$: Since the ``preceding''
path is finished, we start a new one.
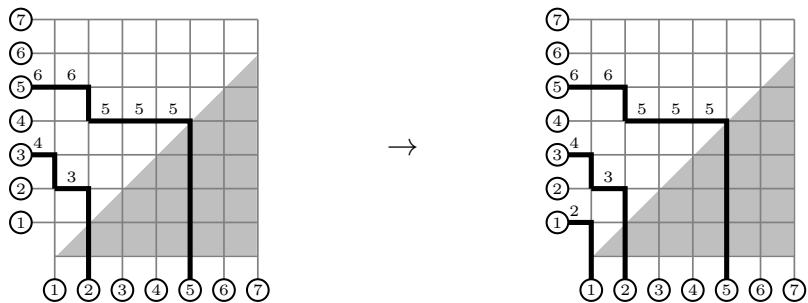
\begin{center}

\hfil
\hbox{\raisebox{2.25cm}{$\to$}}
\hfil
\psset{unit=0.45cm}
\pspicture(-1,-1)(8,9)
\pspolygon[linecolor=white,fillstyle=solid,fillcolor=lightgray](1,1)(7,1)(7,7)
\psset{linewidth=0.05,linecolor=gray}
\multido{\nx=1+1}{7}{\psline(\nx,0)(\nx,8)}
\multido{\ny=2+1}{7}{\psline(0,\ny)(7,\ny)}
\psline(1,1)(7,1)

\psset{fillstyle=none,linecolor=black,linewidth=0.15}
\psline(0,6)(2,6)(2,5)
\rput(0.5,6.35){{\tiny 6}}
\rput(1.5,6.35){{\tiny 6}}
\psline(2,5)(5,5)(5,4)(5,0)
\rput(2.5,5.35){{\tiny 5}}
\rput(3.5,5.35){{\tiny 5}}
\rput(4.5,5.35){{\tiny 5}}
\psline(0,4)(1,4)(1,3)
\rput(0.5,4.35){{\tiny 4}}
\psline(1,3)(2,3)(2,0)
\rput(1.5,3.35){{\tiny 3}}
\psline(0,2)(1,2)(1,0)
\rput(0.5,2.35){{\tiny 2}}
\psset{linewidth=0.06,linecolor=black,fillstyle=solid,fillcolor=white}
\multido{\nx=1+1}{7}{\pscircle(\nx,0){0.35}\rput(\nx,0){{\tiny \nx}}}
\multido{\ny=2+1}{7}{\pscircle(0,\ny){0.35}}
\multido{\ny=1+1}{7}{\rput(0,\ny){\rput(0,1){\tiny\ny}}}

\endpspicture

\end{center}
Reading off the rows from the (heights of the horizontal steps of the)
non--intersecting lattice paths, we obtain the following \DPP\ of dimension $8$
without special parts:
$$
\begin{matrix}
6 & 6 & 5 & 5 & 5 \\
  & 4 & 3 &   & \\
  &   & 2 &   &
\end{matrix}
$$
\end{ex}



\bibliographystyle{plain}
\bibliography{paper}

\end{document}